\documentclass[12pt,a4paper]{amsart}
\usepackage{amssymb}

\usepackage{graphicx,amssymb,amsfonts,epsfig,amsthm,a4,amsmath,url}
\usepackage[latin1]{inputenc}

\newtheorem{thm}{Theorem}
\newtheorem{cor}[thm]{Corollary}
\newtheorem{lem}[thm]{Lemma}

\theoremstyle{definition}

\theoremstyle{remark}
\newtheorem{rem}[thm]{Remark}

\numberwithin{equation}{section}

\newcommand{\M}{\mathcal{M}}

\newcommand{\Supp}{\text{Supp}}

\newcommand{\HH}{\mathcal{H}}
\newcommand{\LL}{\mathcal{L}}
\newcommand{\CC}{\mathcal{C}}
\newcommand{\eps}{\varepsilon}
\newcommand{\SL}{\textnormal{SL}}

\title[Relative Property T for semidirect products]{A characterization of relative Kazhdan Property~T for semidirect products with abelian groups}
\date{January 27, 2010}

\author{Yves de Cornulier}
\address{IRMAR, Campus de Beaulieu, 35042 Rennes CEDEX, France}
\email{yves.decornulier@univ-rennes1.fr}
\author{Romain Tessera}
\address{\'ENS Lyon, 46, all\'ee d'Italie, 69364 Lyon CEDEX~07, France}
\email{Romain.Tessera@umpa.ens-lyon.fr}
\begin{document}


\baselineskip=16pt

\maketitle

\begin{abstract}
Let $A$ be a locally compact abelian group, and $H$ a locally compact group acting on $A$. Let $G=H\ltimes A$ be the semidirect product, assumed $\sigma$-compact.
We prove that the pair $(G,A)$ has Kazhdan's Property~T if and only  if the only countably approximable $H$-invariant mean on the Borel subsets of the Pontryagin dual $\hat{A}$, supported at the neighbourhood of the trivial character, is the Dirac measure. 
\end{abstract}

\section{Introduction}
Let $G$ be a locally compact group and $A$ a subgroup. Recall that the pair $(G,A)$ has {\it Kazhdan's Property~T} (or {\it relative Property~T}, or {\it Property~T}) if every unitary representation of $G$ with almost invariant vectors admits a non-zero $A$-invariant vector. We refer to the book \cite{BHV} for a detailed background.

In this paper, we focus on the special case where $G$ is written as a semidirect product $H\ltimes A$, and $A$ is abelian. Any unitary representation of such a group can be restricted to $A$ and we can then use the spectral theorem to decompose it as an integral of characters. It was thus soon observed that relative Property~T for the pair $(G,A)$ is related to restrictions on invariant probabilities on the Pontryagin dual $\hat{A}$ of $A$. This was first used by D.~Kazhdan \cite{Kaz} in the case of $\SL_n(\mathbf{R})\ltimes\mathbf{R}^n$ for $n\ge 2$. These ideas were then used in a more systematic way, notably by G.~Margulis \cite{Mar} and M.~Burger \cite{Bur}. It was in particular observed that if $H$ is any locally compact group with a representation on a finite-dimensional vector space $V$ over a local field, then $(H\ltimes V,V)$ has Property~T if and only if $H$ does not preserve any probability measure on the Borel subsets of the projective space $\mathbf{P}(V^*)$ over the dual of $V$ (see \cite[Prop.~3.1.9]{corr} for the general statement; the ``if" part follows from \cite[Prop.~7]{Bur}). The idea of using means (i.e.~finitely additive probabilities) instead of probabilities is due to Y.~Shalom \cite[Theorem~5.5]{Sha}, who proved that if $H$ preserves no invariant mean on $\hat{A}-\{1\}$, then $(H\ltimes A,A)$ has Property~T and used related ideas in \cite{Sha2} to prove Property~T for such pairs as $(\SL_2(\mathbf{Z}[X])\ltimes \mathbf{Z}[X]^2,\mathbf{Z}[X]^2)$. Our main result gives the first sufficient condition for relative Property~T in terms of invariant means, which is also necessary.

We say that a Borel mean $m$ on a locally compact space $X$ is {\it countably approximable} if there exists a countable set $\{\nu_n:n\ge 0\}$ of Borel probability measures, whose weak-star closure in $\LL^\infty(X)^*$ (each probability measure being viewed as a mean) contains $m$.

\begin{thm}\label{mainintro}
Let $G=H\ltimes A$ be a $\sigma$-compact locally compact group and assume that the normal subgroup $A$ is abelian. We have equivalences
\begin{itemize}
\item[($\neg$T)]
The pair $(G,A)$ does not have Kazhdan's Property~T.
\item[(M)] There exists a countably approximable $H$-invariant mean $m$ on $\LL^\infty(\hat{A}-\{1\})$ such that $m(V)=1$ for every neighbourhood $V$ of $\{1\}$.
\item[(P)] There exists a net of Borel probability measures $(\mu_i)$ on $\hat{A}$ such that
\begin{itemize}
\item[(P1)] $\mu_i\to\delta_1$ (weak-star convergence in $\CC_c(\hat{A})^*)$;
\item[(P2)] $\mu_i(\{1\})=0$;
\item[(P3)] for every $h\in H$,
$\|h\cdot\mu_i-\mu_i\|\to 0$, uniformly on compact subsets of $H$.
\end{itemize}
\end{itemize}
\end{thm}

Here, Condition (P1) means that $\mu_i(V)\to 1$ for every neighbourhood $V$ of 1 in $\hat{A}$. Also note that since $G$ is assumed $\sigma$-compact, the net in (P) can be replaced by a sequence. In the case of discrete groups, the implication ($\neg$T)$\Rightarrow$(P) has been independently obtained by A.~Ioana \cite[Theorem~6.1]{Ioa}, while its converse was obtained by M.~Burger \cite[Prop.~7]{Bur}.



\begin{cor}\label{coi}
If $H_1\to H$ is a homomorphism with dense image between $\sigma$-compact locally compact groups, then $(H\ltimes A,A)$ has Property~T if and only if $(H_1\ltimes A,A)$ does.

Moreover, if $(H\ltimes A,A)$ has Property~T, then we can find a finitely generated group $\Gamma$ and a homomorphism $\Gamma\to H$ such that $(\Gamma\ltimes A,A)$ has Property~T.
\end{cor}

The first statement of Corollary \ref{coi} shows that, in a strong sense, relative Property~T for such a semidirect product only depends on the image of the action map $H\to\textnormal{Aut}(A)$, and does not detect if this action, for instance, is faithful. It typically applies when $H$ is discrete and $H_1$ is a free group surjecting onto $H_1$. 

\begin{cor}\label{col}
The equivalence between ($\neg$T) and (P) holds for $G$ locally compact (without any $\sigma$-compactness assumption).
\end{cor}

The implication ($\neg$T)$\Rightarrow$(M), which uses standard arguments (similar to \cite[Theorem~5.5]{Sha}), is borrowed, in the discrete case, from \cite[Section~7.6]{Cor}, and improves it in the case when $A$ is not discrete. 

We finish this introduction by giving a relative version of Theorem \ref{mainintro}, generalizing ideas from \cite{CI}. The relative result acutally follows as a corollary from the proof of Theorem \ref{mainintro}. In what follows, all positive functions are assumed to take the value 1 at the unit element. We denote by $\widehat{\mu_i}$ the Fourier-Stieljes transform of $\mu_i$, which is the positive definite function on $A$ defined as $\widehat{\mu_i}(a)=\int\chi(a)d\mu_i(\chi)$ (see \cite[Appendix~D]{BHV}). If $X\subset G$ is a closed subset ($G$ any locally compact group), we also say that $(G,X)$ has relative Property~T if for positive definite functions on $G$, convergence to 1 uniformly on compact subsets of $G$ implies uniform convergence in restriction to $X$ (this extends the previous definition when $X$ is a subgroup, see \cite{corr}). At the opposite, $(G,X)$ has the relative Haagerup Property if there exists positive definite functions on $G$, arbitrary close to 1 for the topology of uniform convergence on compact subsets, but whose restriction to $X$ are $C_0$, i.e.~vanish at infinity. If $G$ is $\sigma$-compact, then this is equivalent to the existence of an affine isometric action on a Hilbert space, whose restriction to $X$ is proper: the proof uses the same argument as the original proof by Akemann and Walter of the equivalence between the unitary and affine definition of Haagerup's Property~\cite{AW}. 

\begin{thm}Under the assumptions of Theorem \ref{mainintro}, assume that $H$ is discrete and suppose that $X\subset A$ is a closed subset. Then we have equivalences 
\begin{itemize}
\item
The pair $(G,X)$ does not have Kazhdan's Property~T.
\item There exists a net of Borel probability measures $(\mu_i)$ on $\hat{A}$ satisfying (P) and such that the convergence of $\widehat{\mu_i}$ to 1 on $X$ is not uniform.
\end{itemize}
If moreover $X$ is $H$-invariant, we have equivalences
\begin{itemize}
\item
The pair $(G,X)$ has relative Haagerup's Property.
\item There exists a net of Borel probability measures $(\mu_i)$ on $\hat{A}$ satisfying (P), with $\widehat{\mu_i}$ is $C_0$ on $X$.
\end{itemize}
\label{relmain}
\end{thm}

In particular, we deduce the following corollary, which generalizes \cite[Theorem~3.1]{CI}.

\begin{cor}
Under the assumptions of Theorem \ref{mainintro}, assume that $H$ is discrete and that $\Lambda$ is a normal subgroup of $H$ whose action on $A$ is trivial. If $X\subset A$, then $(H\ltimes A,X)$ has relative Property~T if and only if $(H/\Lambda\ltimes A,X)$ has relative Property~T. If moreover $X$ is $H$-invariant, then $(H\ltimes A,X)$ satisfies relative Haagerup's Property if and only if $(H/\Lambda\ltimes A,X)$ does.\qed
\end{cor}

\begin{rem}
It is of course better when the condition of Theorem \ref{relmain} on the Fourier-Stieltjes transforms can be made explicit. When $X=A$, wee actually see that, for Borel probability measures on $\hat{A}$, the uniform convergence of $\widehat{\mu_i}$ to one is equivalent to the condition $\mu_i(\{1\})\to 1$. This extends to the case of a subgroup $B$ of $A$ (not necessarily $H$-invariant), by the statement: the convergence of $\widehat{\mu_i}$ to 1 is uniform on $B$ if and only if $\mu_i(H^\perp)\to 1$.

The condition that $\widehat{\mu_i}$ is $C_0$ on $A$ is not easy to characterize but has been long studied (see for instance \cite{Ry}). In view of the Riemann-Lebesgue Lemma, it can be viewed as a weakening of the condition that $\mu_i$ has density with respect to the Lebesgue measure.
\end{rem}

\medskip

To prove the different equivalences, we need to transit through various properties analogous to (P), essentially differing in the way the asymptotic $H$-invariance is stated. Theorem \ref{main} below states all these equivalences and encompasses Theorem \ref{mainintro}. Several of these implications borrow arguments from the proof of the equivalence between various formulations of amenability \cite[Appendix~G]{BHV}. Section \ref{equi} begins introducing some more definitions, notably concerning means, measures, and convolution, and then formulates Theorem \ref{main}. Section~\ref{proofs} contains all proofs.


\section{Equivalent formulations of relative property T for semidirect products}\label{equi}

We need to introduce some notation. 
Let $X=(X,\mathcal{T})$ be a measurable space. Recall that a mean on $X$ is a finitely additive probability measure on the measurable subsets of $X$.
We denote by $\mathcal{L}^{\infty}(X)$ the space of bounded measurable Borel functions on $X$, endowed with the supremum norm $\|\cdot\|_\infty$.
Recall that any mean on $X$ can be interpreted as an element $\bar{m}\in \mathcal{L}^{\infty}(X)^*$ such that $\bar{m}(1)=1$ and $\bar{m}(\phi)\geq 0$ for all non-negative $\phi\in \mathcal{L}^{\infty}(X)$, characterized by the condition $\bar{m}(1_B)=m(B)$ for every Borel subset $B$. By a common abuse of notation, we generally write $m$ instead of $\bar{m}$, and similarly $\mu(f)$ instead of $\int f(x)d\mu(x)$ when $\mu$ is a measure on $X$, and $f$ is an integrable function. Note that any mean $m$ on $X$ can be approximated, in the weak-star topology, by a net $(\nu_i)_{i\in I}$ of finitely supported probabilities (i.e.~finite convex combinations of Dirac measures).

We fix a Haar measure $\lambda$ for $H$. We use the notation $\int f(h)dh$ for the integral of $f\in L^1(H)$ against $\lambda$.
Let $X$ be a measurable space with a measurable action $H\times X\to X$ of $H$.
For every mean $\nu$ on $X$, $h\in H$, and $B$ Borel subset of $X$, we write $(\nu\cdot h)(B)=\nu(hB)$.
Let $\mathcal{UC}_H(X)$ be the subspace of $\mathcal{L}^{\infty}(X)$ whose elements $\phi$ satisfy that $h\to h\cdot \phi$ is continuous from $H$ to $\mathcal{L}^{\infty}(X)$. 
We also need to consider the convolution product $$f* \phi(x)=\int f(h)\phi(h^{-1}x)dh$$ between functions in $L^1(H)$, or between $f\in L^1(H)$ and $\phi \in\LL^{\infty}(X)$. Note that in the first case, $f* \phi\in L^1(H)$, whereas in the second case, $f* \phi\in\LL^\infty(X)$ (see also Lemma \ref{cou}).
If $\mu$ is a measure on $X$, we can define the convolution product 
of $\mu$ and $f\in L^1(H)$ by $\mu*f(B)=\mu(f*1_B)$. Using the Lebesgue monotone convergence theorem, we see it is $\sigma$-additive. It follows (using again the Lebesgue monotone convergence theorem for $\phi\ge 0$) that for all $\phi\in\LL^\infty(X)$ we have
$$(\mu*f)(\phi)=\mu(f*\phi).$$

Let $Y$ be a locally compact Hausdorff space, endowed with its $\sigma$-algebra of Borel subsets. Let $\M(Y)$ be the Banach space of signed Borel regular measures on $Y$ (``regular" is redundant when $Y$ is metrizable), equipped with the total variation norm (i.e.~the norm in $\mathcal{C}_c(Y)^*=\M(Y)$). Note that for $f\in L^1(H)$ and $\mu\in\M(Y)$, we have $\|\mu*f\|\leq \|f\|_1\|\mu\|$.

Let $L^1(H)_{1,+}$ be the subset of $L^1(H)$ consisting of non-negative elements of norm $1$.
Let $\CC_c(H)_{1,+}$ be the set of non-negative, continuous, compactly supported functions $f$ on $H$ such that $\int f(h)dh=1$. Note that $L^1(H)_{1,+}$ and $\CC_c(H)_{1,+}$ are stable under convolution.

\begin{thm}
Let $G=H\ltimes A$ be a $\sigma$-compact locally compact group, with $A$ abelian. Equivalences:
\begin{itemize}
\item[($\neg$T)]  the pair $(G,A)$ does not have Property~T.
\item[(M)] There exists a countably approximable $H$-invariant mean $m$ on $\hat{A}-\{1\}$ such that $m(V)=1$ for every neighbourhood $V$ of $\{1\}$.
\item[(MC)] There exists a Borel $\sigma$-finite measure $\gamma$ on $\hat{A}-\{1\}$ and a mean $m$ on $\LL^\infty(\hat{A}-\{1\})$ belonging to $L^\infty(\hat{A},\gamma)^*$, such that $m(V)=1$ for every neighbourhood $V$ of $\{1\}$, and such that for all $f\in \CC_c(H)_{1,+}$ and $\phi\in \LL^{\infty}(\hat{A})$, 
$$m(f*\phi)=m(\phi).$$
\item[(P)] There exists a net of Borel probability measures $(\mu_i)$ on $\hat{A}$ satisfying (P1), (P2), (P3). 
\item[(PC)] There exists a net of Borel probability measures $(\mu_i)$ on $\hat{A}$ satisfying (P1), (P2), (P3c), where (P3c) is defined as: $\|\mu_i*f-\mu_i\|\to 0,$  for all $f\in \CC_c(H)_{1,+}.$ 
\item[(PQ)] There exists a net of Borel probability measures $(\mu_i)$ on $\hat{A}$ satisfying (P1), (P2), (P3), with the additional property that $\mu_i$ is $H$-quasi-invariant for every $i$.
\end{itemize}
In (P) and (PQ), the net can be chosen to be a sequence. 
Besides, when $G$ is a $\sigma$-compact locally compact group and $A$ a closed abelian normal subgroup (not necessarily part of a semidirect decomposition), then ($\neg$T) implies all other properties (with $H=G/A$), which are equivalent.
\label{main}
\end{thm} 

\begin{rem}
If $(H\ltimes A,A)$ does not have Property~T, we do not necessarily have a net of probabilities $(\mu_i)$, as in any of the properties in Theorem \ref{main}, with {\it density with respect to the Haar measure}. A simple counterexample is given by $\SL_2(\mathbf{R})\ltimes(\mathbf{R}^2\times\mathbf{R})$ (with the trivial action on $\mathbf{R}$), or its discrete analogue $\SL_2(\mathbf{Z})\ltimes (\mathbf{Z}^2\times\mathbf{Z})$. Indeed, we could push this sequence forward to $\mathbf{R}^2$ (resp.~$(\mathbf{R}/\mathbf{Z})^2$) and contradict relative Property~T for $\SL_2(\mathbf{R})\ltimes\mathbf{R}^2$ and $\SL_2(\mathbf{Z})\ltimes\mathbf{Z}^2$.
\end{rem}

\begin{rem}
We could define (M') as the following weak form of (M): there exists an $H$-invariant mean $m$ on $\hat{A}-\{1\}$ such that $m(V)=1$ for every neighbourhood $V$ of $\{1\}$. It can easily be shown to be equivalent to (P'), defined as the existence of a net of Borel probability measures $(\mu_i)$ satisfying (P1),(P2), and (P3'), where (P3') is defined as: $\mu_i-h\mu_i$ tends to zero in the weak-star topology of $\LL^\infty(\hat{A})^*$. We are not able to determine if 
these properties imply ($\neg$T).
\end{rem}

\section{Proof of the results}\label{proofs}
In this section, we first develop a few preliminary lemmas, which hold in a more general context. Then we prove Theorem \ref{main}, and the corollaries.

\begin{lem}Let $X$ be measurable space with measurable action of $H$. For all $f\in L^1(H)$ and for all $\phi\in\LL^{\infty}(X)$, we have $f*\phi\in\mathcal{UC}_H(X)$.\label{cou}\end{lem}
\begin{proof}
If $h\in H$, we have $h\cdot (f*\phi)=(h\cdot f)*\phi$. Therefore, if $h'\in H$ we get
$$\|h\cdot(f*\phi)-h'\cdot(f*\phi)\|_\infty=\|(h\cdot f-h'\cdot f)*\phi\|_\infty$$
$$\le \|h\cdot f-h'\cdot f\|_1\|\phi\|_\infty.$$
Since the left regular action of $H$ on $L^1(H)$ is continuous, we deduce that $g\mapsto g\cdot(f*\phi)$ is continuous from $G$ to $\LL^\infty(X)$, that is, $f*\phi\in\mathcal{UC}_H(X)$.
\end{proof}



\begin{lem} If $A$ is $\sigma$-compact, Condition (P1) is equivalent to:
\begin{itemize}
\item[(P1')] for $a\in A$, we have $\int \chi(a)d\mu_i(\chi)\to 1$, uniformly on compact subsets of $A$.\end{itemize}
\label{conc}
\end{lem}
\begin{proof}
This appears as \cite[Theorem~3.3]{Par} under the assumption that $A$ is second countable (and actually the proof extends to any locally compact abelian group $A$); however we give here a much shorter proof. 

Suppose that (P1) holds. Let $K$ be a compact subset of $A$. There exists a neighbourhood $V$ of $1$ in $\hat{A}$ such that $|1-\chi(a)|\le\eps$ for all $\chi\in V$ and $a\in K$. For $i$ large enough, $\mu_i(V)>1-\eps$, which implies, for all $a\in K$ $$\left|1-\int\chi(a)d\mu_i(\chi)\right|\le \int|1-\chi(a)|d\mu_i(\chi)$$
$$\le  \int_V |1-\chi(a)|d\mu_i(\chi)+\int_{V^c} |1-\chi(a)|d\mu_i(\chi)\le 2\eps.$$

The converse follows from the following claim: for every neighbourhood $V$ of $1$ in $\hat{A}$ and every $\eps>0$, there exists $\eta>0$ and a compact set $K$ in $A$ such that for every Borel measure $\mu$ on $\hat{A}$ satisfying $\sup_{a\in K}|1-\int\chi(a)|d\mu(\chi)\le\eta$, we have $\mu(V)\ge 1-\eps$.

Let us prove this claim. Let $\phi$ be a positive function in $L^1(A)$ with $\int\phi(a)da=1$ (this exists because $A$ is $\sigma$-compact). Set $F(\chi)=\int\phi(a)\chi(a)da$; this is the Fourier transform of $\phi$. In particular, by the Riemann-Lebesgue Lemma, $F$ is continuous and vanishes at infinity. Moreover, $F(1)=1$ and since $\phi>0$, $|F(\chi)|<1$ for all $\chi\neq 1$. Therefore there exists $\rho>0$ such that $\{|F|\ge 1-\rho\}$ is contained in $V$.

Define $\eta=\rho\eps/3$. Let $K$ be a compact neighbourhood of $1$ in $A$ such that $\int_K\phi(a)da\ge 1-\eta$. Let $\mu$ be a Borel probability on $\hat{A}$ such that $$\left|1-\int\chi(a)d\mu(\chi)\right|\le \eta$$
for all $a\in K$. Set $\sigma(a)=\int(1-\chi(a))d\mu(\chi)$.
We have $$\left|\int\phi(a)\sigma(a)\right|\le \left|\int_K\phi(a)\sigma(a)\right|+\left|\int_{K^c}\phi(a)\sigma(a)\right|$$
$$\le \eta+2\eta=3\eta.$$

On the other hand,
$$\int\phi(a)\sigma(a)da=1-\int\left(\int\phi(a)\chi(a)d\mu(\chi)\right)da;$$
since the term in the double integral is summable, we can use Fubini's Theorem, giving
$$\int\phi(a)\sigma(a)=1-\int F(\chi)d\mu(\chi),$$
where 

$$1-\int\phi(a)\sigma(a)=\int F(\chi)d\mu(\chi)$$
$$=\int_{\{|F|> 1-\rho\}} F(\chi)d\mu(\chi)+\int_{\{|F|\le 1-\rho\}} F(\chi)d\mu(\chi),$$
thus
$$\left|1-\int\phi(a)\sigma(a)\right|\le (1-\mu(\{|F|\le 1-\rho\}))+(1-\rho)\mu(\{|F|\le 1-\rho\})$$
$$=1-\rho\mu(\{|F|\le 1-\rho\})$$
so 
$$\left|\int\phi(a)\sigma(a)\right|\ge 1-|1-\int\phi(a)\sigma(a)|\ge
\rho\mu(\{|F|\le 1-\rho\}).$$
Combining with the previous inequality, we obtain.
$$\mu(\{|F|\le 1-\rho\})\le 3\eta/\rho=\eps,$$
hence $$\mu(V)\ge 1-\eps.\qedhere$$
\end{proof}

\begin{lem}\label{Claim1}
Let $X$ be a measurable space with a measurable action of $H$, and $m$ a mean on $\mathcal{UC}_H(X)$. For all $\phi\in \mathcal{UC}_H(X)$ and $f\in L^1(H)$, we have $$
m(f*\phi)=\int f(h)m(h\cdot\phi)dh.$$\end{lem}
\begin{proof}
Fix some $\eps>0$. Let $W$ be a neighbourhood of $1\in H$ such that for every $h\in W,$
\begin{equation}\label{continuityW} 
\|h\cdot \phi-\phi\|_{\infty} \leq \eps.
\end{equation} 
We can write, in $L^1(H)$, $f$ approximately as a finite sum of functions with small disjoint support, namely $f=\sum_{i=1}^kf_i+f_0$ with $\Supp(f_i)\subset h_iW$ for some $h_i\in H$ (when $i\neq 0$) and $\|f_0\|_1\le\eps$ and $\|f\|_1=\sum_j\|f_j\|_1$. Write for short $^h\phi$ for $h\cdot\phi$.

For given $i\neq 0$, we have
\begin{eqnarray*}
& & \left|\int f_i(h)m(^h\phi)dh-m\left(\int f_i(h)\,^h\!\phi\,dh\right)\right|\\
\le &  & \left|\int f_i(h)m(^h\phi)dh-\int f_i(h)m(^{h_i}\phi)dh\right|\\
& + & \left|\int f_i(h)m(^{h_i}\phi)dh-m\left(\int f_i(h)\,^h\!\phi \,dh\right)\right|\\
= &  & \left|\int f_i(h)(m(^h\phi- {^{h_i}\!\phi}))dh\right|+ \left|m\left(\int f_i(h)(^h\phi- {^{h_i}\!\phi})dh\right)\right|\\
\le & & 2\|f_i\|_1\eps\end{eqnarray*}
and
$$\left|\int f_0(h)m(^h\phi)dh-m\left(\int f_0(h)\,^h\!\phi\,dh\right)\right|\le 2\eps\|\phi\|_\infty$$
If we sum over $i$, we deduce
$$\left|\int f(h)m(^h\phi)dh-m\left(\int f(h)\,^h\!\phi\,dh\right)\right|\le 2(\|f\|_1+\|\phi\|_\infty)\eps.$$
Since this holds for any $\eps$, we deduce
$$m(f*\phi)=m\left(\int f(h)\,^h\!\phi\,dh\right)=\int f(h)m(^h\phi)dh.\qedhere$$
\end{proof}

\begin{lem}\label{cl2}
Let $X$ be a measurable space with a measurable action of $H$ by homeomorphisms, and $m$ an $H$-invariant mean on $\mathcal{UC}_H(X)$. Fix $f_0\in  \CC_c(H)_{1,+}$, define a mean by 
$$\tilde{m}(\phi)=m(f_0*\phi), \; \phi\in  \LL^{\infty}(X).$$
Then for all $f\in \CC_c(H)_{1,+}$ and $\phi\in \LL^{\infty}(X)$, 
$$m(f*\phi)=m(\phi).$$
\end{lem}
\begin{proof}
First, $\tilde{m}$ is well-defined by Lemma \ref{cou}. We have to show that $\tilde{m}(f*\phi)=\tilde{m}(\phi)$ for all $f\in  \CC_c(H)_{1,+}$ and $\phi\in \LL^{\infty}(X).$

Let $(f_i)$ be a net in $\CC_c(H)_{1,+}$ with $\Supp(f_i)\to \{1\}$. This implies that 
$\|f *f_i-f\|_1\to 0$, and hence that $\|f*f_i*\phi-f*\phi\|_{\infty}\to 0$, for all $f\in \CC_c(H)_{1,+}$, and $\phi\in \LL^{\infty}(X).$ 
Accordingly $m(f*\phi)= \lim_i  m(f*f_i*\phi)$, which by Lemma \ref{Claim1} equals $\lim_i m(f_i*\phi)$ (since $f_i*\phi\in \mathcal{UC}_H(X)$). This shows that $m(f*\phi)=m(f'*\phi)$ for all $f,f'\in \CC_c(H)_{1,+}$, and all $\phi\in \LL^{\infty}(X)$.
Then for all $f\in \CC_c(H)_{1,+}$ and all $\phi\in  \LL^{\infty}(X),$ $$\tilde{m}(f*\phi)=m(f_0*f*\phi)=m(f_0*\phi)=\tilde{m}(\phi).\qedhere$$
\end{proof}
                                                                                


\begin{proof}[Proof of Theorem \ref{main}]
We are going to prove the implications
$$\text{($\neg$T)}\Rightarrow\text{(P)}\Rightarrow\text{(PQ)}\Rightarrow\text{($\neg$T)}\quad\text{and}$$
$$\text{(P)}\Rightarrow\text{(M)}\Rightarrow\text{(MC)}\Rightarrow\text{(PC)}\Rightarrow\text{(PQ)}\Rightarrow\text{(P)}.$$

\begin{itemize}
\item ($\neg$T)$\Rightarrow$(P).
Let $(\pi,\HH)$ be a unitary representation of $G$ such that                     
$1\prec \pi$ and such that $A$ has no invariant vector.
Let $(K_n)$ be an increasing sequence of compact subsets of $G$ whose interiors cover $G$. Let $(\eps_n)$ be a positive sequence converging to zero.                                                                              
 For each $n$, let $\xi_n$ be a $(K_n,\eps_n)$-invariant                   
vector.                                                                          
Let $E$ be the projection-valued measure associated to                           
$\pi|_{A}$, so that $\pi(a)=\int_{\hat{A}}\chi(a)dE(\chi)$                 
for all $a\in A$. For every $n$, let $\mu_n$ be the                          
probability on $\hat{A}$ defined by $\mu_n(B)=\langle                          
E(B)\xi_n,\xi_n\rangle$.                      
 We have:                                                                         
  $$\|\pi(a)\xi_n-\xi_n\|^2=\int_{\hat{A}}|1-\chi(a)|^2d\mu_n(\chi)\qquad            
\forall a\in A.$$                                                     
Therefore, (P1) results from the almost invariance of $(\xi_n)$.                                                            
Since $\pi$ has no $A$-invariant vector, $\mu_{n}(\{1\})=0$ for               
all $n$. If $f$ is a continuous function on $\hat{A}$, we define a bounded operator $\hat{f}$ on $\mathcal{H}$ by $\hat{f}=\int f(\chi) dE(\chi)$ (actually $\hat{f}$ is the element of the $C^*$-algebra of $\pi|_A$ associated to $f$); note that its operator norm is bounded above by $\|f\|_\infty$. 
For every $h\in H$ and any $f$, we have
$$h\cdot\mu_n(f)=\int fd(h\cdot\mu_n)=\langle \pi(h^{-1})\hat{f}\pi(h)\xi_n,\xi_n\rangle$$
$$=\langle\hat{f}\xi_n,\xi_n\rangle+\langle\hat{f}(\pi(h)\xi_n-\xi_n),\xi_n\rangle$$ $$+\langle\hat{f}\xi_n,\pi(h)\xi_n-\xi_n\rangle+\langle\hat{f}(\pi(h)\xi_n-\xi_n),\pi(h)\xi_n-\xi_n\rangle$$
Thus $$|h\cdot\mu_n(f)-\mu_n(f)|\le 4\|f\|_\infty\|\pi(h)\xi_n-\xi_n\|,$$
so $$\|h\cdot\mu_n-\mu_n\|\le 4\|\pi(h)\xi_n-\xi_n\|$$
which by assumption tends to zero, uniformly on compact subsets of $H$. So (P3) holds.

\item (PQ)$\Rightarrow$($\neg$T). Consider the sequence of Hilbert spaces $\HH_n=L^2(\hat{A}, \mu_n)$, and for every $n$, the unitary action of $H$ on $\HH_n$ defined by 
$$(\pi_n(h)f)(\chi)=f(h\cdot\chi)\left(\frac{d(h\cdot\mu_n)}{d\mu_n}(\chi)\right)^{1/2}.$$ 
There is also a natural action of $A$ on $L^2(\hat{A},\mu_n)$ given by $\pi_n(a)\cdot f(\chi)=\chi(a)f(\chi)$, and since (by a straightforward computation) we have $$\pi_n(h)\pi_n(a)\pi_n(h^{-1})=\pi_n(h\cdot a)\quad \forall h\in H,a\in A,$$so that $\pi_n$ extends to a unitary action of the semidirect product $H\ltimes A$ on $L^2(\hat{A},\mu_n)$. This action has no nonzero $A$-invariant vector. Indeed, let $f$ be an invariant vector. So for every $a\in A$, there exists a Borel subset $\Omega_a\subset\hat{A}$ with $\mu_n(\Omega_a)=1$ and for all $\chi\in\Omega_a$,
$$(\chi(a)-1)f(\chi)=0.$$

If $a\in A$, define its orthogonal $K_a=\{\chi:\chi(a)=1\}$ for all $a\neq 0$. Recall that we assume that $A$ is $\sigma$-compact. If we assume for a moment that $A$ is also second countable, then $A$ is separable; so there exists a sequence $(a_n)$ in $A$ such that $\bigcap_n K_{a_n}=\{1\}$.
If we set $Z=\{f\neq 0\}$, we get $Z\subset K_{a_n}\cup W$, where $W$ is the complement of $\bigcap_n\Omega_{a_n}$. We deduce that $Z\subset W$, which has $\mu_n$-measure zero. So $f=0$ in $L^2(\hat{A},\mu_n)$. If $A$ is only assumed $\sigma$-compact, we proceed as follows: there exists a second-countable open subgroup $B$ of $\hat{A}$ such that $\mu_n(B)>0$ for $n$ large enough (because $\mu_n$ concentrates on $\{1\}$). So we can work in $B$ as we just did in $\hat{A}$ and thus $L^2(\hat{A},\mu_n)$ has no $A$-invariant vector (at least for $n$ large enough).

An immediate calculation gives, for $a\in A$
$$\|1_{\hat{A}}-\pi_n(a)1_{\hat{A}}\|_{L^2(\hat{A},\mu_n)}=2\textnormal{Re}\left(1-\int\chi(a)d\mu_n(\chi)\right),$$
which tends to zero, uniformly on compact subsets of $A$, when $n\to\infty$, by (P1). On the other hand, for every $h\in H$ we have
$$\|1_{\hat{A}}-\pi_n(h)1_{\hat{A}}\|_{L^2(\hat{A},\mu_n)}=
\int_{\hat{A}}\left|1-\left(\frac{d(h\cdot\mu_n)}{d\mu_n}(\chi)\right)^{1/2}\right|^2d\mu_n(\chi),$$
so using the inequality $|1-\sqrt{u}|\le\sqrt{1-|u|}$ for all $u\ge 0$ we get
\begin{eqnarray*}
\|1_{\hat{A}}-\pi_n(h)1_{\hat{A}}\|_{L^2(\hat{A},\mu_n)} &\le &\int_{\hat{A}}\left|1-\frac{d(h\cdot\mu_n)}{d\mu_n}(\chi)\right|d\mu_n(\chi) \\
& = & \|\mu_n-h\cdot\mu_n\|,
\end{eqnarray*}
which tends to zero, uniformly on compact subsets of $H$, when $n\to\infty$, by (P3).
Accordingly, if we consider the representation $\bigoplus\pi_n$, which has no $A$-invariant vector, then the sequence of vectors $(\xi_n)$ obtained by taking $1_{\hat{A}}$ in the $n$th component, is a sequence of almost invariant vectors.


\item(P)$\Rightarrow$(M). View $\mu_n$ as a mean on Borel subsets of $\hat{A}-\{1\}$. Let $m=\lim_\omega\mu_n$ be an accumulation point ($\omega$ some ultrafilter) in the weak-star topology of $\LL^\infty(\hat{A}-\{1\})$.
(P3) immediately implies that $m$ is $H$-invariant. (P1) implies that $\int\chi(a)dm(\chi)=1$ for all $a\in A$. So for every $\eps>0$, we deduce that $m(\{|\chi-1|<\eps\})=1$. In case $A$ is discrete, since those subsets form a prebasis of the topology of $\hat{A}$, we deduce that $m(V)=1$ for every neighbourhood $V$ of $1$ in $\hat{A}$. Hence (M) follows.

When $A$ is not discrete, we need to appeal to Lemma \ref{conc}, which implies that $\mu_n(V)\to 1$ (hence $m(V)=1$) for every Borel neighbourhood $V$ of 1 in $\hat{A}$. 



\item (M)$\Rightarrow$(MC). Let $m$ be an invariant mean as in (M'). Define $\tilde{m}$ as in Lemma \ref{cl2}, which provides the convolution invariance. Clearly, $\tilde{m}(\{1\})=0$. Besides, if $V$ is a closed subset of $\hat{A}$ not containing 1, we see that $f_0*1_V$ is supported by the closed subset $\Supp(f_0)V$, which neither does contain 1. So $\tilde{m}$ is supported at the neighbourhood of $1$. The argument in the proof of Corollary \ref{coi} shows that $\tilde{m}$ also lies in the closure of a countable set $\{\nu_n:n\ge 0\}$ of probability measures 
on $\hat{A}-\{1\}$. If we set $\gamma=\sum 2^{-n}\nu_n$, then $\nu_n$, viewed as a mean, belongs to $L^\infty(\hat{A}-\{0\},\gamma)^*$ (i.e.~vanishes on $\gamma$-null sets), so $m$ also lies in $L^\infty(\hat{A}-\{0\},\gamma)^*$.


 
\item (MC)$\Rightarrow$(PC). Let $m$ be a mean as in (MC) and let ($\nu_i$) be a net of Borel probabilities on $\hat{A}-\{1\}$, converging to $m$ in $\LL^{\infty}(\hat{A})^*$ for the weak-star topology, with $\nu_i$ having density with respect to $\gamma$. We can suppose that $\gamma$ is a probability measure.
Let us show that for any $\eps>0$, and any finite subset $\Omega$ of $\CC_c(H)_{1,+}$, one can find an element $\mu$ in $$V=\left\{\nu\in\M(\hat{A}\setminus\{1\}):\text{Re}\left(\int\chi(a)d\nu(\chi)\right)\ge 1-\eps\right\},$$ such that  $\|\mu*f-\mu\|\leq \eps$ for all $f\in \Omega$. This is exactly, in view of Lemma \ref{conc}, what is required to produce a net $(\mu_i)$ satisfying (PC). First define $\gamma'=\gamma+\sum_{f\in\Omega}\gamma*f$, so that $\gamma'(\{0\})=0$ and each $\mu_i*f$ belongs to $L^1(\hat{A}-\{0\},\gamma')$. 
For every $f\in \CC_c(H)_{1,+}$, the net $(\nu_i*f-\nu_i)$ converges to $0$ for the weak-star topology in $\LL^\infty(\hat{A})^*$, and for every $a\in A$, $\int \chi(a)d\nu_i(\chi)\to 1$.
Since $\gamma'$ is $\sigma$-finite, the dual of $L^1(\hat{A},\gamma)$ is contained in $\LL^\infty(\hat{A})$ (equal to $L^\infty(\hat{A},\gamma')$).
So the convergence of $(\nu_i*f-\nu_i)$ to 0 holds in $L^1(\hat{A}-\{0\},\gamma')$.
Note that $V$ is a closed and convex subset of $\M(\hat{A}\setminus\{1\})$.
Fix $i_0$ such that for all $i\ge i_0$, we have $\nu_i\in V$. Consider the (finite) product 
$$E=L^1(\hat{A}\setminus \{1\},\gamma)^{\Omega},$$ 
equipped with the product of norm topologies. 
Let $\Sigma$ be the convex hull of $$\{ (\nu_{i}*f-\nu_{i})_{f\in \Omega}, \; i\geq i_0\} \subset E.$$
Since $(\nu_i*f-\nu_i)$ converges to $0$ in the weak topology of $E$, the convex set $\Sigma$ contains $0$ in its weak closure. As $E$ is locally convex, by Hahn-Banach's theorem\footnote{The Hahn-Banach Theorem works because we work with the weak topology (and not the weak-star). This is the reason why we need all the measures $\nu_i$ to have density with respect to a given measure $\gamma$. We are not able to bypass this argument.}, the weak closure of $\Sigma$ coincides with its closure in the original topology of $E$. Hence there exists $\mu$ in the convex hull of $\{\nu_i:i\ge i_0\}$ such that $\|\mu*f-\mu\|\leq \eps$ for all $f\in \Omega$; since $V$ is convex, we have 
$\mu\in V$.

\item (PC)$\Rightarrow$(PQ). Let $(\mu_i)$ be as in (PC). By density of compactly supported continuous functions, for all $f\in L^1(H)_{1,+}$, we have $\|f\ast\mu_i-\mu_i\|\to 0$. This convergence is uniform each compact subsets $K$ of $L^1(H)_{1,+}$: this is a trivial consequence of the fact that $(f,\mu)\mapsto\mu$ is 1-Lipschitz for every $\mu$.

Now fix $f_0\in L^1(H)_{1,+}$ and set $\mu'_i=f_0\ast\mu_i$. It is easy to check that it satisfies (P1) and (P2).

By a direct computation, we have, for any $h\in H$ and $\nu\in\M(\hat{A})$, the equality $h\cdot (\nu*f)=\Delta(h)\nu*f_0^h$, where $f_0^h$ is the {\it right} translate of $f_0$, given by $f_0^h(g)=f_0(gh)$. Note that $\Delta(h)f_0^h\in L^1(H)_{1,+}$.
Then for $h\in H$ we have 
$$\|h\cdot\mu'_i-\mu'_i\|=\|h\cdot (\mu_i*f_0)-\mu_i*f_0\|$$
$$=  \|\mu_i*(\Delta(h)f_0^h)-\mu_i*f_0\|$$
$$\le \|\mu_i*(\Delta(h)f_0^h)-\mu_i\|+\|\mu_i*f_0-\mu_i\|.$$
Since the right regular representation of $H$ on $L^1(H)$ is continuous, the function $h\mapsto \Delta(h)f_0^h$ is continuous as well so maps compact subsets of $H$ to compact subsets of $L^1(H)_{1,+}$; therefore the above term converges to zero, uniformly on compact subsets of $H$. So $(\mu'_i)$ satisfies (P3).

Now suppose that we have chosen $f_0>0$ everywhere; this is possible since $H$ is $\sigma$-compact. Let us show that $(\mu'_i)$  satisfies (PQ): it only remains to prove that each $\mu'_i$ is quasi-invariant. Since $h\cdot\mu'_i=\mu_i*(\Delta(h)f_0^h)$, we have to show that the measures $\mu_i*f$, for positive $f\in L^1(H)$, all have the same null sets.
If $B$ is a Borel subset of $\hat{A}$ and $x\in\hat{A}$, we have
\begin{eqnarray*}f*1_B(x)=0 & \Leftrightarrow & \int f(h)1_B(h^{-1}x)dh =0\\ & \Leftrightarrow & \lambda(\{h:f(h)1_B(h^{-1}x) \neq 0\})=0\\ & \Leftrightarrow & \lambda(\{h:1_B(h^{-1}x) \neq 0\})=0\end{eqnarray*}
(since $f$ does not vanish) and this condition does not depend on $f$, provided $f>0$. Thus we have \begin{eqnarray*} \mu_i*f(B)=0 & \Leftrightarrow & \mu_i(f*1_B)=0\\
 & \Leftrightarrow & \mu_i(\{x:f*1_B(x)\neq 0\})=0\end{eqnarray*}
and this condition does not depend on $f$. So $\mu_i*f$ and $\mu_i*f'$ are equivalent measures.
\item (PQ)$\Rightarrow$(P) is trivial.
\end{itemize}

Let us justify the statement about nets and sequences for (P) (the proof for (PQ) being the same). Since $G$ is assumed $\sigma$-compact, there is an increasing sequence $(K_n)$ of compact subsets whose interiors cover $G$. In view of Lemma \ref{conc}, Condition (P) can be written as: for every $\eps>0$ and every $n$, there exists a Borel probability $\mu_{n,\eps}$ on $\hat{A}-\{1\}$ such that $\int \chi(a)d\mu_i(\chi)\ge 1-\eps$ for all $a\in K_n$. So the sequence $(\mu_{n,1/n})$ satisfies the required properties.

For the last statement, first observe that the proof of ($\neg$T)$\Rightarrow (P)$ works without assuming that $A$ is part of a semidirect decomposition. Now all properties except ($\neg$T) only refer to the action on $A$, so their equivalence follows from the theorem applied to the semidirect product $(G\ltimes A,A)$.
\end{proof}

\begin{proof}[Proof of Corollary \ref{coi}]
We use Characterization (M). The ``if" part is trivial. Conversely, suppose that $(H_1\ltimes A,A)$ does not have Property~T. So there exists an $H_1$-invariant mean on $\LL^\infty(\hat{A})$, with $m(1_{\{0\}})=0$ and $m=\lim_\omega\nu_n$ with $\nu_n(V_n)=0$ for some neighbourhood $V_n$ of $1$. Consider the restriction $m'$ of $m$ to $\mathcal{UC}_H(\hat{A})$. Since the action of $H$ on $\mathcal{UC}_H(\hat{A})$ is separately continuous (that is, the orbital maps $H\to \mathcal{UC}_H(\hat{A})$ are continuous), the action on $(\mathcal{UC}_H(\hat{A}),\text{weak*})$ is continuous as well. So the stabilizer of $m'$ is closed in $H$; since it contains the image of $H_1$ in $H$, this shows that $m'$ is $H$-invariant. Fix $f\in C_c(H)_{1,+}$. Thanks to Lemma \ref{cou}, we can define, for $\phi\in\mathcal{L}^\infty(\hat{A})$,
$$m''(\phi)=m'(f*\phi).$$
Clearly, $m''$ is an $H$-invariant mean on $\hat{A}$. Moreover, $m''(1_{\{1\}})=m'(f*1_{\{1\}})=m'(1_{\{1\}})=0$, so $m''$ is not the Dirac measure at 1. Finally we have $m''=\lim_\omega\nu'_n$ in the weak-star topology, where $\nu'_n(\phi)=\nu_n(f*\phi)$, and $\nu'_n$ is a probability on $\hat{A}-\{1\}$. 

For the second statement, assume that $(H\ltimes A,A)$ has Property~T. There exists a compact normal subgroup $K$ in $G=H\ltimes A$ such that $G/K$ is separable \cite[Theorem~3.7]{Com}. Consider a countable subgroup $S$ of $G$ whose image into $G/K$ is dense, and let $T$ be the closure of $S$ in $G$. Set $K'=K/(A\cap K)$. Then $G/K'$ is generated by $T/K'$ and $K/K'$. Since $K/K'$ centralizes $A/K'$, the means preserved by $G$ and by $T$ on the Pontryagin dual of $A/K'$ are the same. This Pontryagin dual is an open subgroup of $\hat{A}$, so the means preserved by $G$ and by $T$ at the neighbourhood of 1 in $\hat{A}$ are the same. Therefore $(T\ltimes A,A)$ has Property~T. Now by the first statement of the Corollary, if $S$ is endowed with the discrete topology, then $(S\ltimes A,A)$ has Property~T. Finally by \cite[Theorem~2.5.2]{Cor}, there exists a finitely generated subgroup $\Gamma$ of $S$ such that $(\Gamma\ltimes A,A)$ has Property~T.
\end{proof}

\begin{proof}[Proof of Corollary \ref{col}]
We first deal with the case when $A$ is not $\sigma$-compact. First, this condition easily implies ($\neg$T) (see for instance \cite[Lemma~2.5.1]{Cor}). It also implies (P). Indeed, let $(G_i)$ be an increasing net of open, $\sigma$-compact subgroups of $G$ and $A_i=G_i\cap A$, $H_i=A_i\cap G_i$. Let $\mu_i$ be the Haar measure on the orthogonal of $A_i$ in $\hat{A}$; note that $\mu_i$ is $H_i$-invariant and $\mu_i(\{1\})=0$ since $A_i$ has infinite index in $A$. So $(\mu_i)$ satisfies (P).

Now suppose that $A$ is $\sigma$-compact. If either ($\neg$T) or (P)
is true for $H\ltimes A$, then it also holds for $L\ltimes A$ for any open subgroup $L$ of $H$. Let us check that conversely, if
it fails for $H\ltimes A$, then it fails for some $\sigma$-compact open subgroup $L\ltimes A$, so that the corollary reduces to the $\sigma$-compact case from the theorem. This is immediate for (P).
For ($\neg$T), if $(H\ltimes A,A)$ has Property~T, by \cite[Theorem~2.5.2]{Cor}, there exists an open compactly generated subgroup $L$ of $H$, containing $A$, such that $(L,A)$ has Property~T.
\end{proof}

\begin{proof}[Proof of Theorem \ref{relmain}]
In either case, suppose that the first condition is satisfied. We have a net $(\varphi_i)$ of positive definite functions on $G$, converging to 1 uniformly on compact subsets of $G$, satisfying some additional condition on $X$. The proof of ($\neg$T)$\Rightarrow$(P) of Theorem \ref{main} constructs a net of Borel measures $(\mu_i)$ on $\hat{A}$, with $\widehat{\mu_i}=\varphi_{|A}$ and $\|\mu_i-h\mu_i\|\to 0$. So we exactly get the second condition.

Conversely, suppose that the second condition is satisfied. Let $\Gamma$ be the subgroup generated by an arbitrary finite subset $S$ of $H$. Denote by $T$ the average operator by $S$. Then $\widehat{T\mu}=T\widehat{\mu}$. If $\widehat{\mu_i}$ is $C_0$ on $X$ and $X$ is $H$-invariant, then $\widehat{T\mu_i}$ is also $C_0$ and $T\mu_i$ is also $\Gamma$-quasi-invariant, we can then follow the proof of (PQ)$\Rightarrow$($\neg$T) of Theorem \ref{main} to obtain a net $(\varphi_i)$ of positive definite functions on $\Gamma\ltimes A$ whose restriction to $A$ is $\widehat{\mu_i}$. 

On the other hand, suppose that the convergence of $\widehat{\mu_i}$ to 1 is not uniform on $X$. Then the convergence of $\widehat{T\mu_i}$ to one is also non-uniform on $X$ (by an obvious positivity argument using that positive definite functions are bounded by one). Again, apply the proof of (PQ)$\Rightarrow$($\neg$T) to obtain the desired net. 

In both cases, we obtain a net on a subgroup of the form $\Gamma\ltimes A$. These functions can be extended to positive definite functions \cite[Exercise~C.6.7]{BHV} on $H\ltimes A$ by taking the value zero elsewhere. If we define the resulting functions as a net indexed by both the indices $i$ and $\Gamma$, the resulting net exactly gives the relative Haagerup Property for $(G,X)$, resp.~the negation of relative Property~T.
\end{proof}


\begin{thebibliography}{KM98b}
\bibitem[AW]{AW} Charles A. {\sc Akemann}, Martin E. {\sc Walter}.
\newblock {\em Unbounded negative definite functions}. \newblock
Canad. J. Math. {\bf 33}, no 4, 862-871, 1981.

\bibitem[BHV]{BHV} B. {\sc Bekka}, P. {\sc de la Harpe}, and A. {\sc Valette}. \newblock  ``Kazhdan's Property~(T)".  \newblock New mathematical monographs: 11, Cambridge University Press, 2008.

\bibitem[Bur]{Bur} M. {\sc Burger}. \newblock {\em Kazhdan constants 
for $\SL_3(\mathbf{Z})$}, J. Reine Angew. Math. {\bf 431}, 36-67, 1991.

\bibitem[CI]{CI} I. {\sc Chifan}, A. {\sc Ioana}. {\em On Relative Property (T) and Haagerup's Property}. ArXiv:0906.5363, to appear in Trans. Amer. Math. Soc.

\bibitem[Com]{Com} W. {\sc Comfort}. \newblock {\em Topological 
groups}. \newblock 1143--1263 in: ``Handbook of Set-Theoretic Topology"
(K. Kunen and J. E. Vaughan, Edt.), North Holland, Amsterdam, 1984.

\bibitem[Cor]{Cor} Y. {\sc Cornulier}. {\em On Haagerup and Kazhdan properties}. Th\`ese EPFL, no 3438 (2006). Dir.: Peter Buser, A. Valette.

\bibitem[Cor2]{corr} Y. {\sc Cornulier}. Relative Kazhdan Property. Annales Sci. \'Ecole Normale Sup. 39(2), 301-333, 2006.

\bibitem[Ioa]{Ioa} A. {\sc Ioana}. {\em Relative Property~(T) for the subequivalence relations induced by the action of $\SL_2(\mathbf{Z})$ on $\mathbf{T}^2$}. Preprint 2009, arXiv 0901.1874 (v1).

\bibitem[Kaz]{Kaz} D. {\sc Kazhdan}. \newblock {\em On the connection 
of the dual space of a group with the structure of its closed 
subgroups}. Funct. Anal. Appl. {\bf 1}, 63-65, 1967.

\bibitem[Mar]{Mar} G. {\sc Margulis}. \newblock ``Discrete
subgroups of semisimple Lie groups''. \newblock Springer, 1991.

\bibitem[Par]{Par} K.~{\sc Parthasarathy}.
``Probability measures on metric spaces". 
Probability and Mathematical Statistics, No. 3 Academic Press, Inc., New York-London, 1967.

\bibitem[Ry]{Ry} R.~Ryan.
Fourier transforms of certain classes of integrable functions. 
Trans. Amer. Math. Soc. {\bf 105}, 1962 102--111. 

\bibitem[Sha]{Sha} Y. {\sc Shalom}. \newblock
{\em Invariant measures for algebraic actions, Zariski dense
subgroups and Kazhdan's property~(T)}. \newblock Trans. Amer.
Math. Soc. \textbf{351}, 3387-3412, 1999.

\bibitem[Sha2]{Sha2} Y. {\sc Shalom}. \newblock  {\em Bounded generation and Kazhdan's Property~(T)}. \newblock Publ. Math. Inst. Hautes \'Etudes Sci. {\bf 90}, 145-168, 1999. 

\end{thebibliography}
\end{document}